\documentclass[11pt]{article}
\usepackage{times,mathptm}
\usepackage{amsfonts,amscd,amssymb}
\input diagrams
\newtheorem{definition}{Definition}[section]
\newtheorem{lemma}[definition]{Lemma}
\newtheorem{proposition}[definition]{Proposition}
\newtheorem{corollary}[definition]{Corollary}
\newtheorem{remark}[definition]{Remark}

\newtheorem{theorem}[definition]{Theorem}

\newcommand{\thlabel}[1]{\label{th:#1}}
\newcommand{\thref}[1]{Theorem~\ref{th:#1}}
\newcommand{\selabel}[1]{\label{se:#1}}
\newcommand{\seref}[1]{Section~\ref{se:#1}}
\newcommand{\lelabel}[1]{\label{le:#1}}
\newcommand{\leref}[1]{Lemma~\ref{le:#1}}
\newcommand{\prlabel}[1]{\label{pr:#1}}
\newcommand{\prref}[1]{Proposition~\ref{pr:#1}}
\newcommand{\colabel}[1]{\label{co:#1}}
\newcommand{\coref}[1]{Corollary~\ref{co:#1}}
\newcommand{\relabel}[1]{\label{re:#1}}

\newcommand{\delabel}[1]{\label{de:#1}}

\newcommand{\eqref}[1]{(\ref{eq:#1})}

\newenvironment{proof}{{\it Proof.}}{\hfill $ \square $ \vskip 4mm}
\newcommand{\Hom}{{\rm Hom}\,}

\newcommand{\Ker}{{\rm Ker}\,}
\newcommand{\im}{{\rm Im}\,}

\newcommand{\End}{{\rm End}\,}

\newcommand{\findim}{{\rm fin.dim}}
\newcommand{\Findim}{{\rm Fin.dim}}
\newcommand{\projdim}{{\rm proj.dim}}

\def\P{{\rm P}}

\def\text#1{\mbox{{\rm #1}}}

\def\dul#1{\underline{\underline{#1}}}

\def\ot{\otimes}

\def\doublerightleft#1#2{{\lower.2ex\vbox{
\hbox{${\smash{\mathop{\longrightarrow}\limits^{#1}}}$}\vspace*{-4mm}
\hbox{${\smash{\mathop{\longleftarrow}\limits_{#2}}}$}}}}

\begin{document}
\title{Separable bimodules and approximation\thanks{Research
supported by the bilateral project BIL99/43 ``New computational, geometric
and algebraic methods applied to quantum groups and
diffferential operators" of the
Flemish and Chinese governments.}}
\author{
S. Caenepeel\\ Faculty of Applied Sciences\\
Free University of Brussels, VUB\\ B-1050 Brussels, Belgium\\
scaenepe@vub.ac.be\and
Bin Zhu\thanks{The second author was partially supported by NSF
10001017 and by the Scientific Foundation for returned overseas Chinese
scholars, Ministry of education. The second author also wishes to
thank the Free University of Brussels for its warm hospitality 
during his visit to Brussels.}\\
Department of Mathematical Sciences\\ Tsinghua University\\
100084 Beijing, China\\ bzhu@math.tsinghua.edu.cn}
\date{}
\maketitle

\begin{abstract}
Using approximations, we give several characterizations of separability
of bimodules. We also discuss how separability properties can be used
to transfer some representation theoretic properties from one ring
to another one: contravariant finiteness of the subcategory
of (finitely generated) left modules with finite
projective dimension, finitistic dimension, finite representation type,
Auslander algebra, tame or wild representation type. 
\end{abstract}

\section{Introduction}\selabel{0}
The notions of approximation and contravariantly finite subcategory were introduced
and studied  by Auslander 
and Smal{\o} \cite{3} in connection with the study of the existence of almost split
sequences in a subcategory.  It turns out that these notions are important in the
study of representation theory of Artin algebras. For example,  Auslander and Reiten
(cf. \cite{1}, \cite{2}) proved that certain contravariantly finite subcategories 
of a
module category are in one-to-one  correspondance to cotilting modules.\\
Auslander and Reiten (\cite{1}, \cite{2}) showed the image of a functor
having a right adjoint is contravariantly finite, we refer to \cite{20}
for a more general result. Now let $R$ and $T$ be rings, and $M$ a
$(T,R)$-bimodule. Then we have a pair of adjoint functors between the categories
of $R$-modules and $T$-modules, and it follows from the Auslander-Reiten result
that the evaluation map $u_M:\ M\otimes_R{}^*M\rightarrow T$ is a right
Im$(F)-$approximation of $T$. This observation enables us to study
separable bimodules and separable extensions
from the point of view of homological finiteness theory. Separable bimodules
have been introduced by Sugano \cite{18}; there has been a revived interest
recently, see for example \cite{4}, \cite{5}, \cite{10} and \cite{11}.\\
In this note, we will apply approximation theory to study ring extensions.
Another aim is to study representation theoretic properties that are shared
by rings connected by a bimodule.\\
Let $A$ be an Artin algebra. If $\P_{\mbox{s}}^{\infty}(A)$, the category
of finitely generated left $A$-modules with finite projective dimension,
is contravariantly finite in the category $A$-mod of finitely generated
left $A$-modules, then the finitistic dimension of $A$ is finite (see
\cite{1}, \cite{2}, \cite{8}). Bass conjectured that the the finitistic
dimension of $A$ is finite, if $A$ is a finite dimensional algebra over a
field $k$ (cf. \cite{8}). The problem is that $\P_{\mbox{s}}^{\infty}(A)$
is not always contravariantly finite in  $A$-mod (see \cite{2} and \cite{8}),
 so it is important to find algebras for which $\P_{\mbox{s}}^{\infty}(A)$
is contravariantly finite. If $T$ is a biseparable extension of $R$, then
the following properties are shared by $T$ and $R$: contravariant finiteness
of the category of finitely generated modules with finite projective
dimension; finitistic and Finitistic dimension. 
In the situation where $R$ and $T$ are Artin algebras, we have that
$R$ is an Auslander algebra if and only if $T$ is an Auslander algebra.
If $R$ and $T$ are Artin algebras connected by a biseparable $(T,R)$-bimodule, 
then $T$ is of finite representation type if and only if $R$ is of finite
representation type; this generalizes a result of Jans \cite{9} and
of Higman \cite{7}. If two finite dimensional algebras $R$ and $T$ over
an algebraically closed field are connected by a biseparable bimodule,
the $T$ is of tame (resp. wild) representation type if and only if
$R$ is of tame (resp. wild) representation type. Some of these results
have been proved in \cite{16}, in the case of skew group ring extensions;
a skew group ring extension is a biseparable extension if the order of
the group is invertible (compare to \cite[Theorem 1.4]{16}).\\
Our paper is organized as follows: in \seref{1}, we recall some preliminary results.
In \seref{2}, we present some characterizations of separability of bimodules,
using approximations, and we discuss how separable bimodules can be used
to construct new separable bimodules. \seref{3} is devoted to the study of
representation theoretic properties of rings connected by a (separable)
bimodule. In \seref{4}, we show that approximations are reflected by
separable functors and we show that the conditional expectation of a Frobenius
extension is an approximation.

\section{Preliminary results}\selabel{1}
Let $R$ and $T$ be rings (associative with unit), and let
$M\in {}_T{\cal M}_R$ be a $(T,R)$-bimodule. Then the right and
left duals
$$M^*=\Hom_R(M,R)~~{\rm and}~~{}^*M={}_T\Hom(M,R)$$
are both $(R,T)$-bimodules; the left and right action are respectively
given by
$$rft(m)=rf(tm)~~{\rm and}~~(m)rgt=((mr)g)t$$
for all $r\in R$, $t\in R$, $m\in M$, $f\in M^*$ and $g\in {}^*M$.\\
A ring extension $R/S$ is a ring homomorphism $i:\ S\to R$. $R$ is
then naturally an $S$-bimodule. $R/S$ is called separable if
the multiplication map $R\ot_S R\to R$ splits as a map of $R$-bimodules.\\
For a ring $T$, we consider the following full subcategories of
the category of left $T$-modules ${}_T{\cal M}$:
\begin{itemize}
\item $T$-mod, consisting of finitely generated left
$T$-modules;
\item ${\cal P}^{\infty}(T)$ consisting of modules with finite
projective dimension;
\item ${\cal P}_S^{\infty}(T)$ consisting of finitely generated left
$T$-modules with finite
projective dimension.
\end{itemize}
Let $k$ be a commutative Artin ring. Recall that a $k$-algebra $A$ is
called an Artin algebra if $A$ is finitely generated as a $k$-module.
An Artin algebra is of finite representation type if there are only 
finitely many isomorphism classes of finitely generated indecomposable
left modules. Now we recall the definition of finitistic 
dimension of an Artin algebra \cite{8}: 
\begin{eqnarray*}
\findim(A)&=&\sup\{\projdim(M)~|~M\in {\cal P}_S^{\infty}(T)\}\\
\Findim(A)&=&\sup\{\projdim(M)~|~M\in {\cal P}^{\infty}(T)\}
\end{eqnarray*}
A conjecture of Bass states that the finitistic dimension of a finite
dimensional algebra over a field $k$ is finite, 
see \cite{2} and \cite{8} for an introduction
and some partial results.\\
Now we recall some definitions from \cite{1} and
\cite{2} that we will need in the sequel. Recall first that a
covariant functor $F:\ {\cal C}\to \dul{\rm Sets}$ is called
finitely generated if and only if there exists an object
$X\in {\cal C}$ and a surjective
natural transformation ${\cal C}(X,\bullet)\to F$. A contravariant
functor is finitely generated if the corresponding covariant
functor ${\cal C}^{\rm op}\to \dul{\rm Sets}$ is finitely generated.

\begin{definition}\delabel{1.1}
Let ${\cal C}$ be a full subcategory of the category ${\cal D}$.\\
(i) ${\cal C}$ is called contravariantly finite in ${\cal D}$
if for all $X\in {\cal D}$, the representable functor
$\Hom_{\cal D}(\bullet, X)$ restricted to ${\cal C}$ is finitely
generated as a functor on ${\cal C}$;\\
(ii) ${\cal C}$ is called covariantly finite in ${\cal D}$
if for all $Y\in {\cal D}$, the representable functor
$\Hom_{\cal D}(Y,\bullet)$ restricted to ${\cal C}$ is finitely
generated as a functor on ${\cal C}$;\\
(iii) ${\cal C}$ is called functorially finite in ${\cal D}$
if ${\cal C}$ is co- and contravariantly finite in ${\cal D}$.
\end{definition}

${\cal C}$ is a contravariantly finite subcatgeory of ${\cal D}$
if and only if the following holds: for each $X\in {\cal D}$,
there exists $X_1\in {\cal C}$ and a morphism $f:\ X_1\to X$
such that $\Hom_{\cal C}(\bullet,f):\ \Hom_{\cal C}(\bullet,X_1)\to
\Hom_{\cal D}(\bullet,X)$ is surjective. This means that every map
$\psi:\ C\to X$, with $C\in {\cal C}$, factors through $f$:
$$\psi=f\circ \varphi:\ C\rTo^{\varphi} X_1\rTo^{f} X$$
for some $\varphi:\ C\to X_1$.
The map $f$ is then called
a right ${\cal C}$-approximation of $X$. Observe that a right
${\cal C}$-approximation is not unique. Left ${\cal C}$-approximations
are defined dually. Pairs of adjoint functors induce approximations:

\begin{lemma}\lelabel{1.1} (\cite{1},\cite{2})
Suppose that $F: {\cal C} \to {\cal D}$ has a right adjoint $G$.
Then $\im(F)$, the full subcategory of ${\cal D}$, consisting of
objects of the form $F(C)$ with $C\in {\cal C}$, is contravariantly
finite in ${\cal D}$. For any $X\in {\cal D}$, the counit map
$\varepsilon_X:\ FG(X)\to X$ is a right $\im(F)$-approximation of $X$.
$\im(G)$ is covariantly finite in ${\cal C}$.
\end{lemma}

\begin{proof}
Let $\eta:\ 1_{\cal C}\to GF$ be the unit of the adjunction. Then
for all $C\in {\cal C}$, we have
$$\varepsilon_{F(C)}\circ F(\eta_C)=1_{F(C)}$$
Consider $f:\ F(C)\to X$ in ${\cal D}$. $\varepsilon$ is a natural transformation,
so we have a commutative diagram
$$\begin{diagram}
FGF(C)&\rTo^{FG(f)}&FG(X)\\
\dTo_{\varepsilon_{F(C)}}&&\dTo^{\varepsilon_X}\\
F(C)&\rTo^{f}&X
\end{diagram}$$
We then compute
$$f=f\circ \varepsilon_{F(C)}\circ F(\eta_C)=\varepsilon_X\circ FG(f)\circ
F(\eta_C)$$
and this is exactly the factorization that we need.
\end{proof}

\leref{1.1} has been generalized in \cite{20}: let
${\cal T}$ be a full subcategory contravariantly finite of ${\cal C}$. 
Then $F({\cal T})$, the full category of ${\cal D}$ consisting of
objects isomorphic to some $F(T)$, with $T\in {\cal T}$, is
a contravariantly finite subcategory of ${\cal D}$.

\section{Separable bimodules}\selabel{2}
The aim of this section is to produce some new separable bimodules 
from given separable bimodules. 
Separable bimodules were introduced by Sugano \cite{18} and studied recently
in \cite{4}, \cite{5}, \cite{10} and \cite{11},
among othres. We recall the definition from \cite{5}. Let
$R$ and $T$ be rings. Given a bimodule ${}_TM_R$, there is a natural
$T$-bimodule homomorphism, 
$$u_M:\ M\otimes_R{}^*M\to T,~` u_M(m\otimes f)= (m)f$$

\begin{definition}\delabel{2.1}
$M$ is called a separable bimodule, or $T$ is called $M$-separable over $R$,
if $u_M$ is a split
$T$-$T$-epimorphism.
\end{definition}

\begin{remark}\relabel{2.1}\rm
It is easy to see that $M$ is separable if and only if there exists $e=\sum
m_i\otimes f_i\in M\otimes_R{}^*M$ such that $u_M(e)=1_T$ and 
$te=et$, for all $t\in T$. $e$ is then called a separable element of
$M$ \cite{10}.
\end{remark}

\begin{definition}\delabel{2.2}
A bimodule $_TM_R$ is called
biseparable if $M$ and $M^*$ are separable and $_TM$, $M_R$ are finitely
generated projective modules.
\end{definition}

Assume that $M_R$ is finitely generated projective. Then
the {\sl evaluation map} $M\to {}^*(M^*)={}_R\Hom(M^*,R)$
is an isomorphism. Identifying $M$ and ${}^*(M^*)={}_R\Hom(M^*,R)$,
we find that
$u_{M^*}:\ M^*\otimes_T M\to R$ is the evaluation map given by
$$u_{M^*}(f\otimes m)=f(m)$$

\begin{definition}\delabel{2.3}
A ring extension $R/S$ is called 
biseparable if ${}_RR_S$ and ${}_SR_R$ are biseparable bimodules.
\end{definition}

To a bimodule ${}_TM_R$, we can associate an adjoint pair of functors
$(F=M\ot_R\bullet ,G={}_T\Hom(M,\bullet)$ between the categories ${}_R{\cal
M}$ and ${}_T{\cal M}$ of respectively left $R$-modules and left $T$-modules.
The same formula defines an adjoint pair of functors between the
categories of bimodules ${}_R{\cal M}_T$ and ${}_T{\cal M}_T$.
Using approximations, we now easily find the following characterizations
of the separability of a bimodule.

\begin{theorem}\thlabel{2.2}
Let ${}_TM_R$ be a bimodule. Then the following statements are equivalent.\\
1) $_TM_R$ is separable, that is, $u_M:\ M\otimes_R {}^*M\rightarrow T$ is a
split $T-T-$epimorphism;\\
2) there exists a split epimorphism of $T$-bimodules
$\phi:\ M\otimes_R{}^*M\rightarrow T$;\\
3) there exists a split epimorphism of $T$-bimodules $\phi:\ M\otimes_R
X\rightarrow T$ for some bimodule ${}_RX_T$;\\
4) $T$ is a direct summand $M\otimes _R {}^*M$ as a $T$-bimodule.
\end{theorem}

\begin{proof}
The implications
$1)\Rightarrow 2)\Rightarrow 3)$ and $1)\Rightarrow 4)\Rightarrow 2)$
are obvious, and we are done if we can show that 3) implies 1).
We have seen in \seref{1} that $\im(F)$ is a contravariantly
generated subcategory of ${}_T{\cal M}_T$, and
$$u_M:\ M\otimes _R{}^*M\rightarrow T,~~
u_M(m\otimes f)=(m)f$$
is a right $\im(F)$-approximation of ${}_TT_T$. This means that
for any $T$-bimodule morphism $\phi:\ M\otimes_RX\to T$, there
exists a $T$-bimodule map $\phi_1:\ M\otimes _RX\rightarrow M\otimes _R{}^*M$ 
such  that $\phi =u_M\circ \phi _1$. If $\phi$ is split, then $u_M$ is
also split, and ${}_TM_R$ is separable.
\end{proof}

Let $A/S$ be a ring extension (in other words,
we have a ring homomorphism $i:\ S\to A$). Then we have two
bimodules ${}_AA_S$ and ${}_SA_A$, and $A/S$ is a separable extension if
and only if ${}_AA_S$ is separable, while $A/S$ is a split extension
if and only if ${}_SA_A$ is separable (see \cite{10}). From \thref{2.2},
we immediately obtain the following result.

\begin{corollary}\colabel{2.3}
A ring extension $A/S$ is separable if and only if $A$ is a direct
summand of $A\ot_S A$ as an $A$-bimodule; $A/S$ is split if and only
if $S$ is a direct summand of $A$ as an $S$-bimodule.
\end{corollary}

Let $X$ be a $T$-bimodule. An element $x\in X$ is called faithful if 
$x\cdot t=0$ implies $t=0$. Denote 
$$_TZ^f(X)_T=\{ x\in X~ |~ x \mbox{ is faithful and } x\cdot t=t\cdot x,  \, 
\forall \  t\in T \ \}$$
\medskip

\begin{lemma}\lelabel{2.2} 
Let $M$ be a $T$-$R$-bimodule. Then $M\otimes _R{}^*M$ contains 
a submodule $N$
which is isomorphic to 
$T$ as a $T$-bimodule if and only if
${}_TZ^f(M\otimes _R {}^*M)_T\not=\Phi$
\end{lemma}

\begin{proof}
Suppose that $N$ is a submodule of ${}_TM\otimes _R{}^*M_T$ and $\phi:
T\rightarrow N$ is a $T$-bimodule isomorphism. It is easy to see that
$e=\phi(1)$ is a casimir element of
$M\otimes _R {}^*M$ \cite{10}. If
$et=0$, then
$\phi(t)=0$, and then 
$t=0$. Therefore $e\in {}_TZ^f(M\otimes _R {}^*M)$.\\
Conversely, let $e\in {}_TZ^f(M\otimes _R {}^*M)$. Then $Te$ is a
$T$-subbimodule of $M\otimes _R{}^*M$.  It is easy
to see that $Te$ is isomorphic to $T$ as a $T$-bimodules. 
\end{proof}

We will now discuss how to produce separable bimodules 
from given separable bimodules.

\begin{theorem}\thlabel{2.4} 
Let ${}_TM_R$ be separable and ${}_RN_S$ a bimodule such that the evaluation
map 
$$u_N^{{}^*M}:\  N\otimes_S {}_R\Hom(N,{}^*M)\to {}^*M,~
u_N^{{}^*M}(m\otimes f)=(m)f$$
is a split $R$-$T$-epimorphism. Then
$_TM\otimes_RN_S$ is separable.
\end{theorem}

\begin{proof}
We have an $S$-$T$-bimodule isomorphism
$${}_T\Hom(M\otimes_R N, T)\cong
{}_R\Hom(N, {}_T\Hom(M,T))$$
It follows that we also have a $T$-$S$-bimodule isomorphism
$$(M\otimes_R N)\otimes_S {}^*(M\otimes_R N)\cong
M\otimes_R (N\otimes_S{}_R\Hom(N, {}^*M)$$
From the fact that $u_N^{{}^*M}$ is a split $R$-$T$-epimorphism,
it follows that ${}^*M$ is an $R$-$T$-direct summand of
$N\otimes_S{}_R\Hom(N,{}^*M)$. Then $T$ is a $T$-$T$-direct summand of
of $(M\otimes_R N)\otimes_S {}^*(M\otimes_R N)$, and \thref{2.2}
tells us that $M\otimes_R N$ is separable.
\end{proof}

As a special case, we have the following consequence.

\begin{corollary}\colabel{2.5}
Let ${}_TM_R$ and ${}_RN_S$ be separable bimodules. Then 
$_TM\otimes _RN_S$ is also a separable bimodule.
\end{corollary}

\begin{proof}
Let ${}_RN_S$ be separable. It is easy to see that $u_N^{{}^*M}:\ 
N\otimes _S {}_R\Hom(N,{}^*M)\to {}^*M,~u_N^{{}^*M}(m\otimes f)=
(m)f$ is a split $R-T$-epimorphism (see also the proof of \leref{3.1}),
and it follows from \thref{2.4} that $_TM\otimes _RN_S$ is separable.
\end{proof}

\begin{theorem}\thlabel{2.6}
Let ${}_TM_R$ be separable and ${}_SX_R$ biseparable. 
Then $\Hom_R(M,X)$ is a separable $T$-$S$-bimodule.
\end{theorem}

\begin{proof}
Set $N=\Hom_R(X,M)$, ${}^*N={}_T\Hom(N,T)$.
We will
prove that $N\otimes_S{}^*N$ contains
$T$ as a $T$-bimodule direct summand, and then, by \thref{2.2}, $N$ is
separable. It is easy to see that 
$X\otimes_R{}_T\Hom(M,T)$ is an $S$-$T$-bimodule.
If we can prove that $N\otimes_S (X\otimes_R{}_T\Hom(M,T))$
contains $T$ as a $T$-bimodule direct summand, then it follows from
\thref{2.2} that $T$ is a direct summand of ${}_TN\otimes _S{}^*N_T$. 
From the biseparability of ${}_SX_R$, it follows that
$$\phi:\ \Hom_R(X,M)\otimes_SX\to M,~~
\phi(f\otimes x)=(x)f$$
is a split epimorphism, hence ${}_TM_R$ is a direct summand of
$\Hom_R(X,M)\otimes_S X$.
Then 
$(\Hom_R(X,M)\otimes_SX)\otimes_R{}_T\Hom(M,T)$ contains
$M\otimes_R\otimes_R{}_T\Hom(M,T)$ as a $T$-bimodule direct summand. Since
${}_TM_R$ is separable, ${}_TT_T$ is a direct summand of 
$M\otimes_R{}_T\Hom(M,T)$, and then it is a direct summand
of 
 $(\Hom_R(X,M)\otimes_SX)\otimes_R{}_T\Hom(M,T)$. 
Since 
$$N\otimes_S(X\otimes_R{}_T\Hom(M,T))\cong 
(\Hom_R(X_R,M_R)\otimes_SX)\otimes_R{}_T\Hom(M,T),$$
${}_TT_T$ is a direct summand of 
$N\otimes_S(X\otimes_R{}_T\Hom(M,T))$. The proof is finished.
\end{proof}

\begin{theorem}\thlabel{2.7}
Let $_TM_R$ be a separable bimodule and $N$ an $T$-$R$-bimodule. Then 
 ${}_T(M\oplus N)_R$ is separable. 
\end{theorem}

\begin{proof}
We have an $R$-$T$-bimodule isomorphism 
$${}^*(M\oplus N)\cong {}^*M\oplus {}^*N$$
Therefore we have $T$-bimodule isomophisms
\begin{eqnarray*}
&&\hspace*{-2cm}
 (M\oplus N)\otimes{}_R{}^*(M\oplus N)\cong 
(M\oplus N)\otimes_R({}^*M\oplus {}^*N)\\
& \cong& (M\otimes _R{}^*M)
\oplus(M\otimes _R{}^*N)
\oplus(N\otimes _R{}^*M)
\oplus(N\otimes _R{}^*N)
\end{eqnarray*}
If $M$ is separable, then $T$ is $T$-bimodule direct summand of
$M\otimes _R{}^*M$, by \thref{2.2}. Then $T$ is also a $T$-bimodule
direct summand of $(M\oplus N)\otimes _R{}^*(M\oplus N)$,
and, again by \thref{2.2}, $_T(M\oplus N)_R$ is separable.
\end{proof}

Let $n$ be a positive integer and $M$ a module. 
The direct sum of $n$ copies of $M$ is denoted by $M^n$. 

\begin{theorem}\thlabel{2.8}
$_TM_R$ is separable if and only if ${}_TM^n_R$ is separable.
\end{theorem}

\begin{proof}
One implication is a direct consequence of \thref{2.7}.
Conversely, assume that ${}_TM^n_R$ is separable. 
${}^*(M^n)\cong ({}^*M)^n$, so we have an 
isomorphism
$$M^n\otimes _R{}^*(M^n) \cong (M\otimes _R{}^*M)^{n^2}$$ 
It follows from the separability of $M^n$ that the map 
$$u_{M^n}:\ (M\otimes _R{}^*M)^{n^2}\to T,~~  
u_{M^n}(m^i\otimes f^i)= \sum
(m^i)f^i$$
is a split epimorphism. Here $m^i\otimes f^i\in M\otimes _R{}^*M$ denotes
the $i$-th component of an element 
$(M\otimes _R{}^*M)^{n^2}$. Let $e=(e^i)$ be a 
separability element of $(M\otimes _R{}^*M)^{n^2}$. It is easy to see that the
sum 
$\sum e^i$ of all entries of $e$ is a separability element of $M$. Therefore 
$M$ is a separable bimodule.
\end{proof}

\section{Representations of rings related by a bimodule}\selabel{3}
There have been various studies of properties shared by rings $R$ and $T$ 
related by a bimodule ${}_TM_R$. A precursor of these
studies is the Higman's Theorem \cite{7}, stating that a finite group has
finite representation type in characteristic $p$ if and only if its Sylow
$p$-subgroup is cyclic. This result appeared later as a Corollary
of Jans' Theorem \cite{9}: for an Artin algebra $R\supseteq T$ in a split
separable extension, $R$ has finite representation type if
and only if $T$ has finite representation type.\\
More results of this type can be found in \cite{5}. In this Section.
we are mainly interested in 
representation theoretic properties shared by rings related by a bimodule,
such as: contravariantly finiteness of the subcategory of modules with finite
projective dimension, Finitistic (or finitistic) dimension and
representation types, Auslander algebras. Some of these were discussed in 
\cite{16}, in the special case of a skew group ring extension.
We will generalize Jans' result to biseparable bimodules. 
We will prove the following result, for finite dimensional algebras
$R$ and $T$ over an algebraically closed field : if there exists a
bimodule ${}_TM_R$, then $T$ is of tame (resp. wild) representation type 
if and only if
$R$ is of tame (resp. wild) representation type.
It would be of interest to have
an Auslander-Reiten theory related to bimodules. We refer to \cite{16} for
some results on skew group ring extensions.\\

For a bimodule $_TM_R$, we have the adjoint pair 
$(F=M\otimes _R\bullet ,G={}_T\Hom(M,\bullet))$ from ${}_R{\cal M}$ to
${}_T{\cal M}$.
For a subcategory ${\cal T}$ of ${}_T{\cal M}$, we denote by 
DS$({\cal T})$ the full
subcategory of ${}_T{\cal M}$ consisting of $T$-modules isomorphic to a direct
summand of $X$ in ${\cal T}$,  by DSIm$(F)$ the full subcategory of
${}_T{\cal M}$ consisting of objects isomorphic to a direct summand of
$F(X)$, where $X\in {}_R{\cal M}$.\\

The following elementary Lemma will be a key tool in our subsequent
results.

\begin{lemma}\lelabel{3.1}
Let $_TM_R$ be a separable bimodule. Then ${}_T{\cal M}={\rm DSIm}(F)$.
\end{lemma}

\begin{proof}
For a left $T$-module $N$, consider the left $R$-module
$${}^*M^N={}_T\Hom(M,N)$$
We have a left $T$-module homomorphism
$$u^N:\ M\otimes _R{}^*M^N\to N,~~u^N(m\otimes f)=(m)f$$
Let $e=\sum m_i\otimes f_i$ be a separability element of $M$,
and consider the map
$$v_N:\ N\to M\otimes _R{}^*M^N,~~v^N(x)=\sum m_j\otimes f_j\cdot x$$
Here $f_j\cdot x\in {}^*M^N$ is defined by the formula
$(m)(f_j\cdot x)=((m)f_j)x$, for all $m\in M$. Now we claim:\\
1) $v^N$ is left $T$-linear. Indeed, for all $t\in T$ and $n\in N$,
we have
\begin{eqnarray*}
&&\hspace{-2cm}
v^N(t\cdot x)=\sum m_j\otimes f_j\cdot (tx)=
[(\sum m_j\otimes f_j)\cdot t]\cdot x\\
&=& t(\sum m_j\otimes f_j)\cdot x
=t(\sum m_j\otimes f_j\cdot x)=tv^N(x)
\end{eqnarray*}
2) $u^N\circ v^N= id_N$: for all $x\in N$, we have
$$(u^N\circ v^N)(x)=\sum(m_j)f_j\cdot x=x$$ 
This means that 
$u^N$ is a split epimorphism of left $T$-modules and $N$ is a 
direct summand of $M\otimes _R{}^*M^N$.
\end{proof}

Let $T/R$ be a ring extension, and consider the adjoint pair
$(F=T\otimes _R\bullet,G)$, where $F$ is the induction functor,
and $G$ is the restriction of scalars functor, between the categories
of left $R$-modules and left $T$-modules. We have a second adjoint
pair $(F'={}_RT\otimes_T\bullet, G'={}_R\Hom(T,\bullet))$ between the
categories of left $T$-modules and left $R$-modules.

\begin{proposition}\prlabel{3.2}
Let $T/R$ be a separable extension with $T$ projective
as a right $R$-module and assume that proj.dim$(_RT)<
\infty$. Then 
 DS$F(\P^{\infty}(R))=\P ^{\infty}(T)$. Moreover, if $\P^{\infty}(R) $ is
contravariantly finite in ${}_R{\cal M}$, then $\P ^{\infty}(T)$ is contravariantly
finite in ${}_T{\cal M}$.  
\end{proposition}

\begin{proof}
Let $X\in {}_R{\cal M}$ with 
projective dimension $m$. Since $T_R$ is projective, proj.dim$(_TT\otimes _RX)\le m$.
Then we have $F(\P ^{\infty}(R))\subseteq \P ^{\infty}(T)$, and then DS$F(\P
^{\infty}(R))\subseteq \P ^{\infty}(T)$ because $\P ^{\infty}(T)$ is closed
under taking direct summands. 
Now let $Y\in \P ^{\infty}(T)$. By the change
of rings Theorem (\cite[section 4.3]{19}), 
$\mbox{proj.dim}(_RY)\le \mbox{proj.dim}(_TY)+  \mbox{proj.dim}(_RT)< \infty$.
${}_TT_R$ is separable, so we can apply \leref{3.1}, and we find that
${}_TY$ is a direct summand of $F(Y)$, proving the first assertion.\\
Assume that $\P^{\infty}(R) $ is contravariantly finite in ${}_R{\cal M}$.
By \cite[Theorem 2.1]{20}, $F(\P^{\infty}$$(R))$ is contravariantly finite in 
${}_T{\cal M}$. We will next show that 
$\mbox{DS}F(\P^{\infty}(R))$ is contravariantly finite in ${}_T{\cal M}$.
Let $Y$ be a left $T$-module, and
$Y_0\rightarrow Y$ be a right
$F(\P^{\infty}(R))$-approximation of $Y$. We verify that $Y_0\rightarrow Y$ is
also a right $\mbox{DS}F(\P^{\infty}(R))$-approximation of $Y$. Take
$Z\in  \mbox{DS}F(\P^{\infty}(R))$, and let a morphism $g:\ Z\to Y$.
$Z\in  \mbox{DS}F(\P^{\infty}(R))$, so there exists 
$Z_1\in F(\P^{\infty}(R))$ and a split $T$-monomorphism $i:\ Z\to Z_1$.
Let $\pi$ be a left inverse of $i$. Then there is a $T$-homomorphism
$h:\ Z_1\to Y_0$, with $g\circ \pi=f\circ h$. It follows that
$g=g\circ \pi\circ i= f\circ h\circ i$, and $g$ factors through $f$.
Then $\mbox{DS}F(\P^{\infty}(R))$ is contravariantly
finite, and therefore $\P ^{\infty}(T)$ is contravariantly finite in 
${}_T{\cal M}$.
\end{proof}

\begin{theorem}\thlabel{3.3}
Let $T/R$ be a biseparable extension. Then
\begin{enumerate}
\item $\P^{\infty}(R) $ is contravariantly
finite in ${}_R{\cal M}$ if and only if $\P^{\infty}(T)$ is contravariantly finite in
${}_T{\cal M}$;
\item $\P_{\mbox{s}}^{\infty}(R)$ is contravariantly
finite in ${}_R{\cal M}$ (resp. in $R$-mod) if and only $\P_{\mbox{s}}^{\infty}(T)$ is
contravariantly finite in ${}_T{\cal M}$ (resp. in $T$-mod);
\item $\mbox{Fin.dim }T=\mbox{Fin.dim }R$;
\item $\mbox{fin.dim }T=\mbox{fin.dim }R$.
\end{enumerate}
\end{theorem}

\begin{proof}
One implication of 1) follows immediately from \prref{3.2}. We prove the converse
direction. First we show that
$$\mbox{DS}F'(\P^{\infty}(T))=\P ^{\infty}(R)$$
It is easy to show that $\mbox{DS}F'(\P^{\infty}(T))\subseteq \P ^{\infty}(R)$.
Conversely, take $Y\in \P ^{\infty}(R)$. $T$ is projective as a right
$R$-module, and $\mbox{proj.dim}(_TT\otimes _RY)<\infty$, so
$\mbox{proj.dim}(_RT\otimes _T(_TT\otimes _RY))<\infty$. Using the
fact that ${}_RT_T$ is separable, we obtain that
$Y$ is an $R$-module direct summand of $F'(_TT\otimes _RY)$, that is,
$Y\in \mbox{DS}F'(\P ^{\infty}(T)$.\\
Now suppose that $\P ^{\infty}$$(T)$ is contravariantly finite.
By \cite[Theorem 2.1]{20}, $F'(\P^{\infty}(T))$ is contravariantly finite
in ${}_R{\cal M}$, and then $\mbox{DS}F'(\P^{\infty}(T))$ is contravariantly finite
(compare to the proof of \prref{3.2}). Consequently
$\P ^{\infty}(R)$ is contravariantly finite in ${}{\cal M}$.\\
The proof of 2) is similar to the proof of 1), because $T$ is finitely generated
and projective as a left and right $R$-module. We omit the details.\\
Now we prove 3). Take 
$X\in \P^{\infty}(T) $ with projective dimension $m$.
$\mbox{proj.dim}({}_RT)=0$, hence $\mbox{proj.dim}({}_RX)\leq m$.
and we claim that we have equality: the projective dimension of
$F(X)$ is smaller than the projective dimension of $X$ as an $R$-module,
and, by \leref{3.1}, $X$ is a $T$-direct summand of $F(X)$.\\
It follows that $m\le \mbox{Fin.dim }R$, and $\mbox{Fin.dim }T\leq
\mbox{Fin.dim }R$. We are done if we can show that
$\mbox{Fin.dim }T\geq
\mbox{Fin.dim }R$. Take 
$Y\in {}_R{\cal M}$ with projective dimension $n$. It is easy to see
that $\mbox{proj.dim}(_TT\otimes _RY)\le n$ and
$\mbox{proj.dim}(_RT\otimes _T(T\otimes _RY))\le n$.
If $\mbox{proj.dim}(_TT\otimes _RY)< n$, then 
$\mbox{proj.dim}(_RT\otimes _T(T\otimes _RY))< n$,
and $Y$, being a direct summand, has projective dimension stricty
smaller than $n$, contradicting the assumption on $Y$. We conclude
that $\mbox{proj.dim}(_TT\otimes _RY)= n$ and
$n\leq \mbox{Fin.dim }T$. This shows that
$\mbox{Fin.dim }T\geq\mbox{Fin.dim }R$.\\
The proof of 4) is similar to the proof of 3), using the fact that
$T$ is finitely generated as a left and right $R$-module.
\end{proof}

\begin{theorem}\thlabel{3.4}
Let Suppose $T/R$ be a biseparable extension of Artin algebras. Then
\begin{enumerate}
\item $\mbox{dom.dim }T =\mbox{dom.dim }R$;
\item $T$ is an Auslander algebra if and only if $R$ is an Auslander algebra.
\end{enumerate}
\end{theorem}

\begin{proof}
1) Let 
$$0\rightarrow T\rightarrow I_0\rightarrow I_1\rightarrow \cdots\rightarrow
I_i\rightarrow \cdots$$
be a minimal injective resolution of $T$. The restriction of scalars
functor $G$ preserves injectives, so this resolution is also an
 injective resolution of $T$ as an $R$-module. 
Now ${}_RR$ is a direct summand of
${}_RT$, so we have an injective resolution of $R$
$$0\rightarrow R\rightarrow I'_0\rightarrow I'_1\rightarrow \cdots\rightarrow
I'_i\rightarrow \cdots$$
with $I'_j$ an $R$-module direct summand of $I_j$, for all $j$.\\
If
$\mbox{dom.dim }T=\infty$, that is, every $I_j$ is projective, then $I'_j$
is projective since $I_j$ is projective in ${}_R{\cal M}$, 
and $\mbox{dom.dim }R=\infty$.\\
If  $\mbox{dom.dim }T=n$, then $I_{n}$ is not projective as a $T$-module, and
the same argument as in the case where $\mbox{dom.dim }T=\infty$ shows
that $\mbox{dom.dim }R\geq n$. If $\mbox{dom.dim }R> n$, then the
$R$-injective resolution of $T$ has the property that $I_n$ is projective
as an $R$-module. This implies that $T\ot_R I_n$ is a projective $T$-module,
and $I_n$ is a projective $T$-module, since it is a $T$-direct summand of
$T\ot_R I_n$. This is a contradiction, so it follows that
$\mbox{dom.dim }R= n$, finishing the proof of part 1).\\
2) Recall that an Artin algeba $T$ is an Auslander algebra if and only if
$\mbox{glob.dim }T\leq 2$ and $\mbox{dom.dim }T\geq 2$. From \cite[Theorem 2.6]{5},
we know that the global dimensions of $R$ and $T$ are equal. Combining this
with part 1), we find 2).
\end{proof}

\begin{remark}\rm
\thref{3.4} has been proved in \cite{16}, in the case of skew group ring
extensions.
\end{remark}

\begin{proposition}\prlabel{3.5}
Let $R$ and $T$ be Artin algebras, and assume that there exists a
biseparable $(T,R)$-bimodule $M$. Then $R$ is of finite representation
type if and only if $T$ is of finite representation type.
\end{proposition}

\begin{proof}
Assume that $R$ is of finite representation type. 
It follows from the separability of $M$ and \leref{3.1} that the full
subcategory $\mbox{ind}({}_T{\cal M})$ of 
${}_T{\cal M}$ consisting of finitely generated indecomposable modules 
coincides with the full subcategory $\mbox{ind}(\mbox{DSIm}(F))$ of 
$\mbox{DSIm}(F)$ consisting of finitely generated indecomposable modules.
$\mbox{ind}(\mbox{DSIm}(F))$ is of finite representation type since 
$R$ is of finite representation type, and it follows that $\mbox{ind}(T)$
is of finite representation type. The converse implication follows from
the fact that the $(R,T)$-bimodule $M^*$ is biseparable.
\end{proof}

\begin{corollary}\colabel{3.6}
Let $A$ be a finite dimensional algebra over a field $k$. If there is a separable
bimodule ${}_AM_k $, 
 then $A$ is of finite representation type.
\end{corollary}

\begin{proof}
We note that $\mbox{ind}({}_k{\cal M})$ contains just one object. It then
follows from the proof of \prref{3.5} that
$A$ is of finite representation type.
\end{proof}

From now on we assume that $R$ and $T$ are finite dimensional Artin algebras
over an algebraically closed field $k$. Any finite dimensional $k$-algebra
is Morita equivalent to a basic algebra of the form $kQ/I$.
$Q$ is called the Gabriel quiver of $A$ (cf. \cite{17}).
Recall that a quiver $Q$ is a pair $(Q_0,Q_1)$, where $Q_0$ is the set of vertices
and $Q_1$ is the set of arrows. A
$k$-algebra $A=kQ/I$ is of tame representation type if for each dimension
vector 
$\textbf{\underline{z}}\in \textbf{N}^{Q_0}$, there exist finitely many parametrizing
$A{\hbox{-}}k[t]$-bimodules 
$M_1,\cdots ,M_s$ satisfying the two following conditions:
\begin{enumerate}
\item every $M_i$ is finitely generated and
free as a right $k[t]-$module;
\item every indecomposable $A$-module $X$ for which
$\underline{\mbox{dim}}X=\textbf{\underline{z}}$ is isomorphic to a module of the
form $M_i\otimes (k[t]/(t-\lambda ))$, with $i\in \{ 1,\cdots , s\}$ and
$\lambda
\in k$.
\end{enumerate}
 It was proved in \cite{12},\cite{13} that $A$ is of tame representation type 
if and only if $A$ is weakly tame. This means that for every
$\textbf{\underline{z}}\in
\textbf{N}^{Q_0},$ there is a family of finitely generated 
$A{\hbox{-}}k[t]$-bimodules
$M_1,\cdots ,M_s$ such that each indecomposable
$A$-module $X$
 with $\underline{\mbox{dim}}X=\textbf{\underline{z}}$ is
a direct summand of $M_i\otimes_{k[t]}S$ for some 
$i\in \{ 1,\cdots , s\}$ and a simple $k[t]-$module $S$.

\begin{proposition}\prlabel{3.7}
Let $R$ and $T$ be finite dimensional algebras over an algebraically
closed field.Then 
$R$ is of tame representation type if and only if
$T$ is of tame representation type.
\end{proposition}

\begin{proof}
We divide the proof into two parts. First, we show that we can restrict
attention to the
situation where $T$ and $R$ are basic. Then we prove the Theorem for
basic algebras $R$ and $T$.\\
Let $T'$ be the basic algebra of
$T$. Then there is a $(T',T)$-bimodule $X$ that induces an equivalence $F_1=X\otimes
_T\bullet :\ {}_T{\cal M} \rightarrow {}_{T'}{\cal M}$. We claim that $_{T'}X\otimes
_TM_R$ is biseparable. Firstly, since ${}_{T'}X_T$ and ${}_TM_R$ are separable, by
\coref{2.5} we have that $F_1(M)$ is separable. Secondly, we have an isomorphism
of
$(R,S)$-bimodules:
$$\Hom_R(X\otimes _TM, R)\cong \Hom_T(X, \Hom(M, R))$$ 
${}_{T'}X_T$ is biseparable, because it induces a Morita equivalence.
It follows
from \thref{2.6} that\\
$\Hom_T(X,\Hom_R(M, R))$ is separable. Third,
we can easily verify that $F_1(M)_R$ and ${}_{T'}F_1(M)$ are finitely generated
projective. Then we have that 
$_{T'}X\otimes _TM_R$ is biseparable. Dually, let $R'$ be the basic algebra of $R$
and ${}_RY_{R'}$
the bimodule inducing an equivalence $G_1=\bullet\otimes _RY_{R'}:\
{\cal M}_R\rightarrow{\cal M}_{R'}$. A similar argument shows that 
$F_1(M)\otimes _RY$
is a  biseparable $(T',R')$-bimodule. So without loss of generality, we
can assume that $T$ and $R$ are basic algebras and ${}_TM_R$ is biseparable.\\

Assume that $R=kQ/I$ is of tame representation type. We will prove 
that $T=k\Gamma /J$ is
weakly tame. Let 
$\textbf{\underline{w}}=(w(i))_{i\in \Gamma_0}$ be a dimension vector. 
We prove that there are only finitely many
dimension vectors 
$\textbf{\underline{v}}=(v(j))_{j\in Q_0}$ with 
$\textbf{\underline{v}}=\underline{\mbox{dim}}G(N)$ for some $T$-module $N$ with
$\underline{\mbox{dim}}N=\textbf{\underline{w}}$. Let $P(j)$ be the indecomposable
projective 
$R$-module corresponding to the vertex $j$.
Then we have an isomorphism 
$$ \Hom_R(P(j),G(N))\cong \Hom_T(M\otimes_RP(j), N)$$
$\Hom_T(M\otimes_RP(j), N)$
 is a direct summand of $\mbox{Hom}_T(M, N)$, and it follows that
$$\mbox{dim}_k\Hom_R(P(j),G(N))\leq
\mbox{dim}_k\Hom_T(M,N)$$
for all $j$. The left hand side has an upper bound
$\mbox{dim}k\Hom_k(M, k^{\sum {w(i)}})$. Therefore there are only finitely
many dimension vectors 
$\textbf{\underline{v}}=(v(j))_{j\in Q_0}$ with 
$\textbf{\underline{v}}=\underline{\mbox{dim}}G(N)$ for some $T$-module $N$ with
 $\underline{\mbox{dim}}N=\textbf{\underline{w}}$. Let
${\rm ind}_R(\textbf{\underline{v}})$ denote the subcategory of 
${}_R{\cal M}$ consisting of
indecomposable modules with dimension vector $\textbf{\underline{v}}$.  Since $R$ is
tame, there are finitely many $(R,k[t])$-bimodules
$M_i$ which are finitely generated free right $k[t]]$-modules and parametrize all
${\rm ind}_R(\textbf{\underline{v}})$, where
$\textbf{\underline{z}}$ fullfils the above estimation. Let
 $M_i'=M\otimes _RM_i$. Then we have that any indecomposable left $T$-module $N$ 
with
 $\underline{\mbox{dim}}N=\textbf{\underline{z}}$ is isomorphic to a direct summand 
of $M_i'\otimes_{k[t]}S$ for some simple
$k[t]$-module $S$. This proves that $T$ is weakly tame. 
For the converse direction, we
use the separability of the $(R,T)$-bimodule
 $M^*$ and  the tameness of $T$. The arguments are then the duals of the
ones presented above. 
\end{proof}

Combining Propositions \ref{pr:3.5} and \ref{pr:3.7} with Drozd's tame-wild
dichotomy theorem \cite{6}, we obtain the following result.

\begin{theorem}\thlabel{3.8}
Let $R$ and $T$ be finite dimensional algebras over
an algebraically closed field, and $_TM_R$ a biseparable bimodule.
Then
\begin{enumerate}
\item $R$ is of finite representation type if and only if
$T$ is of finite representation type;
\item $R$ is of tame and infinite representation type if and only if
$T$ is of tame and infinite representation type;
\item $R$ is of wild representation type if and only if
$T$ is of wild representation type.
\end{enumerate}
\end{theorem}

\section{Right approximations}\selabel{4}
It is a difficult problem to decide which subcategories are
contravariantly finite, and to find a right approximation
(see \cite{1}, \cite{2}, \cite{20}). In this Section, we will
see that separable functors reflect approximations.
Then we give some descriptions of split extensions and Frobenius extensions.\\

Let $F:\ {\cal C}\to {\cal D}$ be a covariant functor.
$F$ induces a natural transformation
$${\cal F}:\ \Hom_{\cal C}(\bullet,\bullet)\to 
\Hom_{\cal D}(F(\bullet),F(\bullet)),~{\cal F}_{C,C'}(f)=F(f)$$
Recall from \cite{14} that $F$ is called separable if
${\cal F}$ splits as a natural transformation, that is, there exists
a natural transformation
$${\cal P}:\ \Hom_{\cal D}(F(\bullet),F(\bullet))\to \Hom_{\cal C}(\bullet,\bullet)$$
such that ${\cal P}\circ {\cal F}$ is the identity natural transformation
on $\Hom_{\cal C}(\bullet,\bullet)$. For a detailed study of separable functors,
we refer the reader to \cite{4}.

\begin{proposition}\prlabel{4.1}
Let $F:\ {\cal C}\to {\cal D}$ be a separable functor, 
and ${\cal T}$ a full subcategory of ${\cal C}$. Let $C_1\in {\cal T}$
and $C\in {\cal C}$, and a morphism $f:\ C_1\to C$. If $F(f):\ F(C_1)
\to F(C)$ is a right (resp. a left) $F({\cal T})$-approximation of $F(C)$,
then $f:\ C_1\to C$ is a right (resp. a left) ${\cal T}$-approximation of $C$.
\end{proposition}

\begin{proof}
Assume that 
$F(f):\ F(C_1)\to F(C)$ is a right $F({\cal T})$-approximation of $F(C)$.
Let $g:\ B\to C$ be a morphism in ${\cal C}$. Then there exists a 
morphism $h:\
F(B)\to F(C)$ in ${\cal D}$ such that the following diagram commutes:
$$\begin{diagram}
F(B)&\rTo^{h}&F(C_1)\\
\dTo^{F(1_B)}&&\dTo^{F(f)}\\
F(B)&\rTo^{F(g)}&F(C)
\end{diagram}$$
that is, $F(g)=F(f)\circ h$.
It follows from the separability of $F$ that we have the following commutative
diagram in
${\cal C}$:
$$\begin{diagram}
B&\rTo^{{\cal P}(h)}&C_1\\
\dTo^{1_B}&&\dTo^{f}\\
B&\rTo^{g}& C
\end{diagram}$$
or $g=f\circ {\cal P}(h)$. This means that $f:\ C_1\to C$ is a right
${\cal T}$-approximation of $C$. The proof in the case of a left
approximation is similar.
\end{proof}

We will now study ring extensions from the point of view of
approximations. Consider a ring extension $R/S$. We use the following
notation for the restriction of scalars functors:
$$G':\ {\cal M}_R\to {\cal M}_S~~{\rm and}~~
G:\ {}_S{\cal M}_R\to {}_S{\cal M}_S$$
It is well-known that $G$ and $G'$ have a right adjoint and a left adjoint,
and it follows that $\im(G)$ (resp. $\im(G')$) is covariantly and
contravariantly finite im ${}_S{\cal M}_S$ (resp. ${\cal M}_S$)
(see \seref{1}), and we can construct left and right approximations.
In particular, a right $\im(G)$-approximation of ${}_SS_S$ is the map
$\phi:\ R^*\to S$, mapping $f\in R^*$ to $f(1)$. $\phi$ is also a right 
$\im(G')$-approximation of $S_S$

\begin{proposition}\prlabel{4.4}
If $\phi:\ {}_SR_R\rightarrow \mbox{Hom}_S(R, S)$ is an isomorphism of
$(S,R)$-bimodules, then 
 $E=\phi(1):\ {}_SR_S\to S$
is a right $\im(G)$-approximation of $_SS_S$, and also a right
$\im(G')$-approximation of $S_S$.
\end{proposition}

\begin{proof}
We know that $u_S:\ {}_SR^*\otimes_RR_S\to S$ is a right 
$\im(G)$-approximation of
$_SS_S$. We have the following commutative diagram of $S$-bimodules:
$$\begin{diagram}
R&\rTo^{E}& S\\
\dTo^{\cong}&&\dTo^{=}\\
R\ot_R R&&S\\
\dTo^{\phi\ot 1_R}&&\dTo^{=}\\
R^*\ot_R R&\rTo^{u_S}&S
\end{diagram}$$
All the vertical maps are isomorphisms, so $E$ is a right 
$\im(G)$-approximation of $S$. $u_S$ is also a right 
$\im(G')$-approximation of $S$, and the same is true for $E$.
\end{proof}

\begin{corollary}\colabel{4.5}
If $R/S$ is a Frobenius extension, with Frobenius system $\{ E, x_i,y_i\}$,
then $E:\ R\to S$ is a right $\im(G)$-approximation
of $_SS_S$.
\end{corollary}

A right $\im(G')$-approximation $E:\ R\rightarrow S$ of $S_S$ is called
non-degenerate if $\Ker(E)$
contains no non-zero right ideal of $R$.

\begin{theorem}\thlabel{4.6}
Let the $S$-bimodule map $E:\ R\to S$ be a right 
$\im(G')$-approximation of $S_S$ and $A=\End_S(R)$. Then
$$R^*=\Hom_S(R,S)= E\circ A$$
as $(S,R)$-bimodules.
If $E$ is non-degenerate, then $R\cong E\circ A$
as $(S,R)$-bimodules. Conversely, if $R\cong E\circ A$
as $(S,R)$-bimodules, then there exists a non-degenerate
$\im(G')$-approximation $E_1:\ {}_SR_S\to S$ of $S_S$.
\end{theorem}

\begin{proof}
Let $E:\ R\to S_S$ be a right $\im(G')$-appriximation of $S_S$.
For any $f\in \Hom_S(R,S)$, 
there is an $h:\ R_S\to R_S$ such that $f=E\circ h$, and it follows that
$R^*= E\circ A$, as right $A$-modules. Observe that $E\circ A$ is a left
$S$-module, since $s\cdot (E\circ f)= E\circ (s\cdot f)$, for all
$s\in S$ and $f\in A$. Indeed, for all $x\in R$, we have
$$s\cdot(E\circ f)(x)=s\cdot E(f(x))=E(s f(x))=(E\circ(s\cdot f))(x)$$
Finally, we have a monomorphism $R\to A$, mapping $r$ to $m_r:\ R\to R$,
$m_r(x)=rx$, and we have an $(S,R)$-bimodule isomorphism $R^*=E\circ A$.\\
Now suppose that $E$ is non-degenerate. It is easy to see that
$E\circ R$ is an $(S,R)$-subbimodule of $R^*$. The map $\alpha:\
R\to E\circ R$, $\alpha(x)=E\circ m_r$ is surjective. It is injective,
since $E$ is non-degenerate: if $x\in \Ker \alpha$, then $E(xr)=0$,
for all $r\in R$, and $xR$ is a non-zero right ideal of $R$ contained
in $\Ker E$. Thus $\alpha$ is an isomorphism of right $(S,R)$-bimodules.\\
Conversely, let $\alpha:\ R\to E\circ R$ be an isomorphism of
$(S,R)$-bimodules, with inverse $\alpha^{-1}$. Let $\alpha(1)=E\circ a$
and $\alpha^{-1}(E)=b$. Then 
$$E=\alpha(\alpha^{-1}(E))=\alpha(b)=\alpha(1)b=E\circ m_a\circ m_b
=E\circ m_{ab}$$
This implies that $E(abx)=E(x)$, and $E((ab-1)x)=0$, for all $x\in R$.
Then $ab=1$, otherwise $R(ab-1)$ is a non-zero right ideal contained in
$\Ker E$. Also
$$1=\alpha^{-1}(\alpha(1))=\alpha^{-1}(E\circ m_a)=\alpha^{-1}(E)a=ba$$
$E\circ m_a$ is an $S$-bimodule map since
$$(E\circ m_a)(sx)=\alpha(1)(sx)=(\alpha(1)\cdot s)(x)=\alpha(s)(x)=
s\alpha(1)(x)=s(E\circ m_a)(x)$$
for all $x\in R$ and $s\in S$. $E=E\circ m_{ab}=E\circ m_a\circ m_b$
is a right $\im(G')$-approximation of $S_S$, so $E\circ m_a$ is also
a right $\im(G')$-approximation of $S_S$. Finally, $E\circ m_a$ is
non-degenerate: assume that there exists $x\in R$ such that
 $(E\circ m_a)(xr)=0$ for all $r\in R$. Then 
$$\alpha(x)(r)=(\alpha(1)x)(r)=(E\circ m_a)(xr)=0$$
for all $r\in R$, hence $\alpha(x)= 0$, and $x=0$, since $\alpha$ is
bijective.
\end{proof}

\end{document}